\documentclass{amsart}

\makeatletter
\@namedef{subjclassname@2020}{%
  \textup{2020} Mathematics Subject Classification}
 
\makeatother % Cambiar el ano de la MSC, de 2010 a 2020

\usepackage{amsthm,amssymb,amsfonts,latexsym,mathtools,thmtools,amsmath,lscape}
\usepackage[T1]{fontenc}
\usepackage{tikz-cd} 
\usepackage{enumitem}
\usepackage{longtable}
\usepackage{multirow}
\usepackage{caption}
\usepackage{mathrsfs} 
\usepackage{hyperref} 
\usepackage{float}
\usepackage[all]{xy}
\usepackage{booktabs}
\usepackage{tabularx}
\usepackage{array} % necesario para redefinir columnas

\hypersetup{
    colorlinks=true,
    linkcolor=blue,
    filecolor=blue,      
    urlcolor=blue,
    linktocpage=true
}

\newtheorem{theorem}{Theorem}[section]
\newtheorem{lemma}[theorem]{Lemma}
\newtheorem{proposition}[theorem]{Proposition}

\theoremstyle{definition}
\newtheorem{definition}[theorem]{Definition}
\newtheorem{corollary}[theorem]{Corollary}
\newtheorem{remark}[theorem]{Remark}

\newtheorem{example}[theorem]{Example}

\theoremstyle{remark}

\numberwithin{equation}{section}

\begin{document}
\raggedbottom

\title[Differential geometry of generalized down-up algebras]{Noncommutative differential geometry \\ of generalized down-up algebras}

\author{Andr\'es Rubiano}
\address{Universidad Distrital Francisco Jos\'e de Caldas}
\curraddr{Campus Universitario}
\email{aarubianos@udistrital.edu.co - ORCID: 0009-0009-1633-8018} 

\author{Armando Reyes}
\address{Universidad Nacional de Colombia - Sede Bogot\'a}
\curraddr{Campus Universitario}
\email{mareyesv@unal.edu.co - ORCID: 0000-0002-5774-0822}

\thanks{The second author was supported by Faculty of Science, Universidad Nacional de Colombia - Sede Bogot\'a, Colombia [grant number 65488].}



\date{}

\dedicatory{Dedicated to the memory of Georgia Benkart}

\begin{abstract} 

We investigate the differential smoothness of down-up algebras and generalized down-up algebras over fields of characteristic zero. We construct diagonal calculi for suitable constant-polynomial presentations and triangular calculi for polynomials satisfying two affine eigenvalue identities. The latter construction proves differential smoothness of split homogeneous down-up algebras when the quotient of the roots is $1$ or has infinite multiplicative order. We also give a connected integrable three-dimensional calculus on the universal enveloping algebra of the Lie algebra $\mathfrak{sl}_2$ and determine the limitations of diagonal commutation rules in the two presentations.

\end{abstract}

\maketitle
\enlargethispage{2pt}

%\tableofcontents 

\section{Introduction}

In noncommutative differential geometry, algebras, modules, and differential operators encode geometric data. In the classical setting, functions on a space form an algebra, sections of a vector bundle form a module, and differential forms carry a graded product and an exterior derivative. These algebraic features can be retained when the algebra of functions is replaced by an associative algebra. A choice of differential calculus then provides the setting in which one can formulate connections, curvature and integration (for more details on the subject, see Connes \cite{Connes1994} and Brzezi\'nski \cite{Brzezinski2008}).

For a differential graded algebra $(\Omega A,d)$ and a left $A$-module $M$, a {\em connection} is a linear map $\nabla^0:M\to\Omega^1 A\otimes_A M$ satisfying
\[
\nabla^0(am)=da\otimes_A m+a\nabla^0(m),\quad {\rm for\ all}\ a\in A,\ m\in M.
\]
This formulation preserves the product rule for covariant differentiation without referring to points of an underlying space, while the use of tensor products suggests a complementary construction involving spaces of module homomorphisms. Brzezi\'nski \cite{Brzezinski2008} developed this second approach through the notion of {\em hom-connection}, which acts on suitable dual modules and satisfies a corresponding right Leibniz rule. In the case of the algebra itself, such an operator is usually called a {\em divergence}. This terminology is consistent with the ordinary divergence when one takes the usual differential calculus on a Euclidean space.

The relationship between connections and hom-connections extends beyond a formal analogy. Brzezi\'nski \cite{Brzezinski2008} constructed hom-connections from connections on bimodules and from appropriate connections in noncommutative principal bundles. He also described their extension to higher-degree dual forms and introduced a {\em curvature operator}. When this curvature vanishes, the extended operators define a {\em chain complex of integral forms}. Thus a differential calculus can be accompanied by two related structures: its complex of differential forms and a complex built from their right module duals. {\em Twisted multi-derivations} provide explicit constructions of calculi and divergences in the free-module setting \cite{BrzezinskiElKaoutitLomp2010}. The corresponding treatment of finitely generated projective modules was also developed by him \cite{Brzezinski2011}, allowing the module-theoretic analogue of a cotangent bundle to be nontrivial.

On the other hand, several notions of {\em smoothness} coexist in noncommutative algebra, reflecting different aspects of classical differential geometry. {\em Formal smoothness} can be expressed through a lifting property, or equivalently through projectivity of the bimodule of universal one-forms; in the terminology of Cuntz and Quillen \cite{CuntzQuillen1995}, the resulting algebras are called {\em quasi-free}. {\em Homological smoothness} instead requires a finite resolution of the algebra by finitely generated projective bimodules; this condition underlies the duality framework of Van den Bergh \cite{VandenBergh1998}. These properties are defined through extension and resolution data. The existence of a differential calculus with an appropriate duality is a separate geometric question, and its verification requires explicit control of forms, their multiplication and their differential.

The notion considered in this paper is the {\em differential smoothness} introduced by Brzezi\'nski and Sitarz \cite{BrzezinskiSitarz2017}. This is motivated by the compatibility between differential and integral forms on a smooth orientable manifold. In their algebraic formulation an integrable differential calculus possesses a volume form, for which multiplication of complementary degree forms produces pairings that identify each module of forms with the appropriate right dual. Equivalently, the complex of differential forms is isomorphic - with a possible twist of the left module structure - to a complex of integral forms with the grading reversed, yielding an algebraic analogue of {\em Hodge-star duality}. In this setting, an affine algebra of integer Gelfand--Kirillov dimension $n$ is differentially smooth when such a calculus can be chosen connected and $n$-dimensional. This definition relates the growth of the algebra to the number of available differential degrees.

This framework includes both familiar quantum spaces and deformations of classically singular spaces. Brzezi\'nski and Sitarz \cite{BrzezinskiSitarz2017} treated the noncommutative pillow, quantum cones and quantum lens spaces, as well as the differential structures that connect these examples. Related constructions include the integral calculi of \cite{BrzezinskiElKaoutitLomp2010}, the affine Hopf algebras studied in \cite{Brzezinski2015}, and the generalized Weyl algebras considered in \cite{Brzezinski2016}. Criteria for passing differential smoothness to suitable skew polynomial extensions were established by Brzezi\'nski and Lomp \cite{BrzezinskiLomp2018}. Further works concern three-dimensional skew polynomial algebras, skew PBW extensions and diffusion algebras on an arbitrary finite number of generators \cite{AlmulhemBrzezinski2019, ReyesSarmiento2022, RubianoReyes2025DSSPBWKt,RubianoReyes2025DSSPBWNote,RubianoReyes2026RUMA}.

For bi-quadratic algebras with PBW basis defined by Bavula \cite{Bavula2023}, the diagonal commutations in \cite[Theorem 3.2]{RubianoReyes2026BQAThree} give a way to organize the verification of the differential and the volume form. Related constructions using skew polynomial and generalized Weyl presentations appear in \cite{RubianoReyes2026BQA}. In the {\em down-up setting} introduced by Benkart and Roby \cite{BenkartRoby1998}, the polynomial appearing in the commutator of the up and down generators leads to mixed terms in the action on one-forms. These terms will be used explicitly below.

Benkart and Roby \cite{BenkartRoby1998} introduced {\em down-up algebras} to study relations between operators that move down and up in partially ordered sets. Their defining relations depend on three scalar parameters and are cubic in two generators. With this presentation, this family of algebras contains the enveloping algebras of the three-dimensional Heisenberg Lie algebra and of the Lie algebra $\mathfrak{sl}_2$. They proved that every down-up algebra has a PBW basis and Gelfand--Kirillov dimension three. Kirkman et al. \cite{KirkmanMussonPassman1999} subsequently characterized the Noetherian members of the family: they are precisely those for which the parameter multiplying the reversed cubic term is nonzero, and these are also precisely the domains.

Later, Cassidy and Shelton \cite{CassidyShelton2004} enlarged this class by introducing {\em generalized down-up algebras}. These algebras have three generators, two relations governing the interaction with a distinguished generator, and a third relation involving an arbitrary polynomial in that generator. For an affine polynomial with nonzero linear coefficient, eliminating the distinguished generator gives an ordinary down-up presentation. Allowing higher-degree polynomials includes algebras related to the enveloping algebra of $\mathfrak{sl}_2$, such as those studied by Smith \cite{Smith1990}. Cassidy and Shelton established ring-theoretic and representation-theoretic properties of this larger class, while Carvalho and Lopes \cite{CarvalhoLopes2009} developed its structural description further and determined automorphism groups under conformality and parameter restrictions. Their formulation also supplies the corrected parameter correspondence between generalized and ordinary down-up algebras.

These presentations make the two families natural candidates for a study based on differential calculi. In the Noetherian case, a generalized down-up algebra admits descriptions both as an iterated skew polynomial ring and as a generalized Weyl algebra over a polynomial ring in two variables \cite{CassidyShelton2004,CarvalhoLopes2009}. Such descriptions expose automorphisms and skew derivations that can enter the construction of bimodules of one-forms. Almulhem and Brzezi\'nski \cite{AlmulhemBrzezinski2019} constructed a class of skew derivations on complex Noetherian generalized down-up algebras, providing a direct antecedent for this differential-geometric question. More recently, the limitations of certain ambiskew smoothness criteria for down-up algebras, and the problem of studying generalized down-up algebras were discussed in \cite[Sections 4.5 and 5]{RubianoReyes2026Ambiskew}. These observations motivate a treatment that keeps track of the defining relations, the choice of calculus and the integrability conditions simultaneously.

In this paper, we investigate the differential smoothness for down-up algebras and generalized down-up algebras within this framework. Section \ref{preliminaries} recalls the differential and ring-theoretic preliminaries. In Section \ref{Differentialandintegralcalculusdownup}, we first examine the diagonal construction and then give an Ore extension of the algebra of forms that permits mixed commutations. Theorem \ref{gdu-constant} treats constant polynomials, Theorem \ref{triangular-gdu} gives the affine eigenvalue criterion, and Theorem \ref{TheoremDU} applies it to ordinary down-up algebras. Section \ref{ExamplesDU} presents examples in both families, including a separate construction for $U(\mathfrak{sl}_2)$. Section \ref{Futurework} gives an obstruction for affine diagonal rules and discusses the scope of the results.

Throughout the paper, $\mathbb{N}$ denotes the set of natural numbers including zero. The word ring means an associative ring with identity not necessarily commutative. $Z(R)$ denotes the center of the ring $R$. All vector spaces and algebras (always associative and with unit) are over a fixed field $\Bbbk$. $\Bbbk^*$ denotes the non-zero elements of $\Bbbk$. As usual, the symbols $\mathbb{R}$ and $\mathbb{C}$ denote the fields of real and complex numbers, respectively.

\section{Preliminaries}\label{preliminaries}

We first recall the differential and integral structures used in the definition of differential smoothness. We then introduce down-up algebras and generalized down-up algebras, specifying the conventions needed to pass between their presentations.

\subsection{Differential smoothness of algebras}\label{Difsmoothness}

The material in this subsection follows \cite[Section 2]{BrzezinskiSitarz2017}, together with the theory of hom-connections in \cite{Brzezinski2008,BrzezinskiElKaoutitLomp2010}.

\begin{definition}[{\cite[Section 2.1]{BrzezinskiSitarz2017}}]\label{calculus}
\begin{enumerate}[label=\textup{(\roman*)},leftmargin=*]
\item A \emph{differential graded algebra} is a nonnegatively graded algebra $\Omega=\bigoplus_{k\in\mathbb{N}}\Omega^k$, with multiplication denoted by $\wedge$, equipped with a degree-one $\Bbbk$-linear map $d$ such that $d^2=0$ and
\[
d(\eta\wedge\xi)=d\eta\wedge\xi+(-1)^k\eta\wedge d\xi, \quad {\rm for\ all} \ \eta\in\Omega^k,\ \xi\in\Omega.
\]
\item A \emph{differential calculus over $A$} is a differential graded algebra $(\Omega A,d)$ with $\Omega^0 A=A$ for which
\[
\Omega^k A=\text{span}_{\Bbbk}\{a_0\,da_1\wedge\cdots\wedge da_k\mid a_0,\ldots,a_k\in A\}, 
\]
for every $k\geq 1$. This requirement is called the \emph{density condition}. The calculus is \emph{connected} if $\ker(d|_A)=\Bbbk$.
\item A calculus has \emph{dimension $n$} when $\Omega^n A\ne0$ and $\Omega^k A=0$ for $k>n$. It \emph{admits a volume form} if there is an element $\omega\in\Omega^n A$ that is a free generator of $\Omega^n A$ both as a right and as a left $A$-module. Such an element $\omega$ is called a \emph{volume form}.
\end{enumerate}
\end{definition}

Repeated use of the Leibniz rule moves coefficients in the density condition to the right. Consequently, $\Omega^k A$ is also spanned by the products $da_1\wedge\cdots\wedge da_k\,a_0$. A volume form allows every top-degree form to be expressed uniquely as $\omega a$, and hence determines a right $A$-module isomorphism
\begin{equation}\label{pi-volume}
\pi_{\omega}:\Omega^n A\longrightarrow A,\  \pi_{\omega}(\omega a)=a.
\end{equation}
The left action of $A$ on this generator is encoded by
\begin{equation}\label{nu-volume}
a\omega=\omega\nu_{\omega}(a),\  a\in A.
\end{equation}
Since the same $\omega$ freely generates the top forms on both sides, $\nu_{\omega}$ is an algebra automorphism. If only a free right generator is specified, the same formula defines an endomorphism, which need not be invertible.

For each $k\in\mathbb{N}$, let
\[
\mathcal{I}_k A:=\text{Hom}_A(\Omega^k A,A)
\]
be the right $A$-module dual of $\Omega^k A$. It is an $A$-bimodule under
\begin{equation}\label{dual-actions}
(a\cdot\phi\cdot b)(\eta)=a\phi(b\eta),
\quad {\rm for\ all} \ a, b\in A,\  \phi\in\mathcal{I}_k A,\ \eta\in\Omega^k A.
\end{equation}
The direct sum $\mathcal{I} A=\bigoplus_{k\in\mathbb{N}}\mathcal{I}_k A$ has a right action of the graded algebra $\Omega A$, defined by evaluation after multiplication:
\begin{equation}\label{form-action}
(\phi\cdot\eta)(\xi)=\phi(\eta\wedge\xi),
\quad {\rm for\ all} \ \phi\in\mathcal{I}_{k+m}A, \  \eta\in\Omega^k A, \  \xi\in\Omega^m A.
\end{equation}
This action lowers the degree of the dual form by the degree of $\eta$.

\begin{definition}[{\cite[Definition 2.1]{Brzezinski2008}}]
A \emph{divergence}, or \emph{hom-connection on $A$}, is a $\Bbbk$-linear map $\nabla:\mathcal{I}_1 A\to A$ such that
\begin{equation}\label{divergence}
\nabla(\phi\cdot a) = \nabla(\phi)a+\phi(da),
\quad {\rm for\ all} \ \phi\in\mathcal{I}_1 A, \ a\in A.
\end{equation}
\end{definition}

A divergence determines maps $\nabla_k:\mathcal{I}_{k+1} A\to\mathcal{I}_k A$ through
\begin{equation}\label{extended-divergence}
\nabla_k(\phi)(\eta)=\nabla(\phi\cdot\eta)+(-1)^{k+1}\phi(d\eta),
\quad {\rm for\ all} \ \eta\in\Omega^k A.
\end{equation}
We identify $\mathcal{I}_0 A=\text{Hom}_A(A,A)$ with $A$ by evaluation at $1$; with this identification, $\nabla_0=\nabla$. The extensions obey the rule 
\begin{equation}\label{extended-leibniz}
\nabla_k(\phi\cdot\eta)=\nabla_{k+m}(\phi)\cdot\eta+(-1)^{k+m}\phi\cdot d\eta,
\end{equation}
where $\phi\in\mathcal{I}_{k+m+1} A$ and $\eta\in\Omega^m A$ \cite[Lemma 3.2]{Brzezinski2008}. Taking $k=m=0$ gives \eqref{divergence}.

The divergence is \emph{flat} if its curvature $\nabla\circ\nabla_1:\mathcal{I}_2 A\to A$ vanishes. Under this hypothesis, $\nabla_k\circ\nabla_{k+1}=0$ for every $k\geq0$, and one obtains the \emph{complex of integral forms}
\[
\cdots\longrightarrow\mathcal{I}_2 A\xrightarrow{\nabla_1}\mathcal{I}_1 A\xrightarrow{\nabla}A.
\]
The corresponding \emph{integral} is the canonical quotient map
\[
\Lambda:A\longrightarrow\text{Coker}\nabla=A/\text{Im}\nabla.
\]
This means that flatness is the condition that turns the sequence of extended divergences into a chain complex \cite[Section 3]{Brzezinski2008}.

If $X$ is a left $A$-module and $\nu$ is an automorphism of $A$, the notation ${}^{\nu}X$ refers to the same vector space with left action $a\cdot_{\nu}x=\nu(a)\cdot x$. When $X$ is a bimodule, its right action is left unchanged. This twist is needed to compare differential forms with their right duals in the noncommutative setting.

\begin{definition}[{\cite[Definition 2.1]{BrzezinskiSitarz2017}}]\label{integrable}
An $n$-dimensional differential calculus $(\Omega A,d)$ is \emph{integrable} if it admits a flat divergence and there are an automorphism $\nu$ of $A$ and $A$-bimodule isomorphisms
\[
\Theta_k:\Omega^k A\longrightarrow{}^{\nu}\mathcal{I}_{n-k}A, \quad {\rm for\ all} \ 0\leq k\leq n,
\]
which intertwine the differentials. Explicitly, for these values of $k$, every square
\[
\begin{tikzcd}[column sep=large,row sep=large]
\Omega^k A \arrow[r,"d"] \arrow[d,"\Theta_k"'] &
\Omega^{k+1} A \arrow[d,"\Theta_{k+1}"]\\
{}^{\nu}\mathcal{I}_{n-k}A \arrow[r,"\nabla_{n-k-1}"'] &
{}^{\nu}\mathcal{I}_{n-k-1}A
\end{tikzcd}
\]
commutes, where $\Omega^0 A=A$ and $\mathcal{I}_0 A=A$. The element $\omega=\Theta_n^{-1}(1)$ is an \emph{integrating volume form}.
\end{definition}

The following characterization allows integrability to be checked directly from multiplication of forms.

\begin{proposition}[{\cite[Theorem 2.2]{BrzezinskiSitarz2017}}]\label{integrability}
Let $(\Omega A,d)$ be an $n$-dimensional differential calculus. The following conditions are equivalent.
\begin{enumerate}[label=\textup{(\arabic*)},leftmargin=*]
\item The calculus is integrable.
\item For some algebra automorphism $\nu$, there are $A$-bimodule isomorphisms \linebreak $\Theta_k:\Omega^k A\to{}^{\nu}\mathcal{I}_{n-k}A$, for $0\leq k\leq n$, satisfying
\[
\Theta_{k+m}(\eta\wedge\xi)=(-1)^{(n-1)m}\Theta_k(\eta)\cdot\xi, 
\]
for every $\eta\in\Omega^k A$, $\xi\in\Omega^m A$ and $k+m\leq n$.
\item There exist an algebra automorphism $\nu$ and an $A$-bimodule map $\vartheta:\Omega^n A\to{}^{\nu}A$ such that every map
\[
\ell_{\vartheta}^{k}:\Omega^k A\longrightarrow\mathcal{I}_{n-k}A, \quad \ell_{\vartheta}^{k}(\eta)(\xi)=\vartheta(\eta\wedge\xi),
\]
is bijective for $0\leq k\leq n$.
\item The calculus admits a volume form $\omega$ for which
\[
\ell_{\pi_{\omega}}^{k}:\Omega^k A\longrightarrow\mathcal{I}_{n-k}A, \quad \ell_{\pi_{\omega}}^{k}(\eta)(\xi)=\pi_{\omega}(\eta\wedge\xi),
\]
is bijective for every $1\leq k\leq n-1$.
\end{enumerate}
\end{proposition}

For a volume form, the maps in degrees $0$ and $n$ are already bijective. Hence condition (4) tests only the intermediate degrees. A volume form is {\em integrating} exactly when it satisfies this condition \cite[Remark 2.3]{BrzezinskiSitarz2017}.

Finite projectivity of the modules of forms is useful both geometrically and computationally. The following criteria express it through finite reconstruction formulas, and the second criterion also establishes integrability.

\begin{proposition}[{\cite[Lemmas 2.6 and 2.7]{BrzezinskiSitarz2017}}]\label{dual-bases}
\begin{enumerate}[label=\textup{(\arabic*)},leftmargin=*]
\item Suppose $(\Omega A,d)$ is an integrable $n$-dimensional calculus with integrating volume form $\omega$. Fix $k$. If there are finitely many $\eta_i\in\Omega^k A$ and $\xi_i\in\Omega^{n-k} A$ such that
\[
\eta=\sum_i\eta_i\pi_{\omega}(\xi_i\wedge\eta), \  \eta\in\Omega^k A,
\]
then $\Omega^k A$ is a finitely generated projective right $A$-module.
\item Suppose an $n$-dimensional calculus admits a volume form $\omega$. Assume that, for every $1\leq k\leq n-1$, there are finite families $\omega_i^j,\overline{\omega}_i^j\in\Omega^j A$, for $j=k,n-k$, giving the two expansions
\begin{equation}\label{two-reconstructions}
\eta=\sum_i\omega_i^k\pi_{\omega}(\overline{\omega}_i^{n-k}\wedge\eta) =\sum_i\nu_{\omega}^{-1}\!\left(\pi_{\omega}(\eta\wedge\omega_i^{n-k})\right) \overline{\omega}_i^k
\end{equation}
for all $\eta\in\Omega^k A$. Then $\omega$ is an integrating volume form, and every $\Omega^k A$ is finitely generated and projective as both a left and a right $A$-module.
\end{enumerate}
\end{proposition}

To compare the length of a calculus with the algebra it describes, one uses Gelfand--Kirillov dimension. If $A$ is affine and $V$ is a finite-dimensional generating subspace containing $1$, then
\[
\operatorname{GKdim} A=\limsup_{m\to\infty}\frac{\log\dim_{\Bbbk}V^m}{\log m}.
\]
As is well-known, this value is independent of the choice of $V$ \cite{KrauseLenagan2000}.

\begin{definition}[{\cite[Definition 2.4]{BrzezinskiSitarz2017}}]\label{smoothness}
An affine algebra $A$ with $\operatorname{GKdim} A=n\in\mathbb{N}$ is \emph{differentially smooth} if it has a connected integrable differential calculus of dimension $n$.
\end{definition}

\begin{example}\label{polynomial-calculus}
For $A=\Bbbk[x_1,\ldots,x_n]$, the ordinary polynomial de Rham calculus has
\[
\Omega^k A=\bigoplus_{1\leq i_1<\cdots<i_k\leq n}dx_{i_1}\wedge\cdots\wedge dx_{i_k}\,A, \ \omega=dx_1\wedge\cdots\wedge dx_n.
\]
Products of complementary basis forms show that the maps $\ell_{\pi_{\omega}}^k$ are bijective. In characteristic zero, a polynomial whose partial derivatives all vanish is constant, so the calculus is connected. Since $\operatorname{GKdim} A=n$, the polynomial algebra is differentially smooth. Noncommutative examples over $\mathbb{C}$ include the quantum plane, the standard Podle\'s sphere and the coordinate algebra of $SU_q(2)$, with the usual real deformation parameter $0<q<1$ (e.g.  \cite{BrzezinskiElKaoutitLomp2010,BrzezinskiSitarz2017}).
\end{example}

\begin{remark}[{\cite[Example 2.5]{BrzezinskiSitarz2017}}]
The existence of an appropriate calculus can fail even for a commutative affine algebra. The algebra $\mathbb{C}[x,y]/\langle xy\rangle$ has Gelfand--Kirillov dimension one, but admits no connected integrable one-dimensional differential calculus, so it is not differentially smooth.
\end{remark}

We now turn to the algebraic objects of primary interest in this paper.

\subsection{Down-up algebras}\label{down-up}

As we said in the Introduction, down-up algebras were introduced by Benkart and Roby \cite{BenkartRoby1998}. Let us see the details.

\begin{definition}[{\cite[Section 1]{BenkartRoby1998}; \cite[Section 1]{KirkmanMussonPassman1999}}]\label{down-up1}
Given $\alpha,\beta,\gamma\in\Bbbk$, the \emph{down-up algebra} $A(\alpha,\beta,\gamma)$ is the unital associative $\Bbbk$-algebra with generators $d,u$ and defining relations
\begin{align}
d^2u&=\alpha dud+\beta ud^2+\gamma d,\label{du-first}\\
du^2&=\alpha udu+\beta u^2d+\gamma u.\label{du-second}
\end{align}
Equivalently, it is the quotient of the free associative algebra $\Bbbk\langle d,u\rangle$ by the two-sided ideal generated by the differences of the two sides of these relations.
\end{definition}

No nonvanishing restriction is imposed on the parameters in this definition. Although there are only two algebra generators, the ordered normal forms involve three independent exponents.

\begin{proposition}[{\cite[Theorem 3.1 and Corollary 3.2]{BenkartRoby1998}}]\label{du-pbw}
For every $\alpha,\beta,\gamma\in\Bbbk$, the set
\[
\left\{u^i(du)^j d^k\mid i,j,k\in\mathbb{N}\right\}
\]
is a $\Bbbk$-basis of $A(\alpha,\beta,\gamma)$. Moreover, $
\operatorname{GKdim} A(\alpha,\beta,\gamma)=3$. 
\end{proposition}

With respect to the filtration by word length in $d,u$, the basis elements of filtration degree at most $m$ are those satisfying $i+2j+k\leq m$. If $\gamma=0$, the defining relations are homogeneous for $\deg d=\deg u=1$. These observations explain why a calculus witnessing differential smoothness must have dimension three, even though the presentation in Definition \ref{down-up1} uses two generators.

\begin{proposition}[{\cite[Main Theorem]{KirkmanMussonPassman1999}}]\label{du-noetherian}
For $A = A(\alpha,\beta,\gamma)$, the following conditions are equivalent:
\begin{enumerate}[label=\textup{(\arabic*)},leftmargin=*]
\item $\beta\ne0$;
\item $A$ is left and right Noetherian;
\item $A$ is a domain;
\item the subalgebra $\Bbbk[ud,du]$ is a polynomial algebra in two variables.
\end{enumerate}
When these conditions hold, $A$ is Auslander-regular and has global dimension three.
\end{proposition}

The growth assertion in Proposition \ref{du-pbw} also holds when $\beta=0$. Thus the parameter restriction in Proposition \ref{du-noetherian} concerns ring-theoretic properties rather than the dimension appearing in Definition \ref{smoothness}.

\begin{example}[{\cite[Section 1]{BenkartRoby1998}}]\label{Heisenberg}
Let $\mathfrak h_3$ be the {\em Lie algebra} with basis $x,y,z$ and brackets $[x,y]=z$, $[x,z]=[y,z]=0$. Then $A(2,-1,0)\cong U(\mathfrak h_3)$. 
Indeed, for $c=du-ud$, the two defining relations are given by $[d,c]=0$ and $[c,u]=0$. The correspondence $d\mapsto x$, $u\mapsto y$, and $c\mapsto z$ identifies both presentations. Note that its inverse is induced by $x\mapsto d$, $y\mapsto u$ and $z\mapsto du-ud$.
\end{example}

\begin{example}[{\cite[Section 1]{BenkartRoby1998}}]\label{sl-two}
Let $e,f,H$ be the standard basis of $\mathfrak{sl}_2(\Bbbk)$. Then 
\[
[e,f]=H, \  [H,e]=2e \quad {\rm and} \quad  [H,f]=-2f.
\]
Then $A(2,-1,-2)\cong U(\mathfrak{sl}_2(\Bbbk))$ 
under $d\mapsto e$ and $u\mapsto f$. In fact,
\[
e^2f-2efe+fe^2 = [e,H] = -2e \quad {\rm and} \quad ef^2 - 2fef + f^2e = [H,f] = -2f.
\]
Conversely, in $A(2,-1,-2)$ the elements $d,u,du-ud$ satisfy the displayed Lie brackets. This yields the inverse homomorphism and fixes the sign of the parameter $\gamma$ for this choice of generators.
\end{example}

\subsection{Generalized down-up algebras}\label{generalized-down-up}

Cassidy and Shelton \cite{CassidyShelton2004} introduced generalized down-up algebras by allowing a polynomial relation involving a third generator. We follow the convention used by Carvalho and Lopes \cite{CarvalhoLopes2009}. To distinguish the scalar in this presentation from the parameter $\gamma$ of a down-up algebra, we denote it by $\varepsilon$.

\begin{definition}[{\cite[Section 1]{CassidyShelton2004}; \cite[Section 1.1]{CarvalhoLopes2009}}]\label{generalized-down-up1}
Let $f\in\Bbbk[x]$ and $r,s,\varepsilon\in\Bbbk$. The \emph{generalized down-up algebra} $L(f,r,s,\varepsilon)$ is the unital associative $\Bbbk$-algebra generated by $d,u,h$, subject to
\begin{align}
dh&=rhd-\varepsilon d,\label{gdu-dh}\\
hu&=ruh-\varepsilon u,\label{gdu-hu}\\
du&=sud-f(h).\label{gdu-du}
\end{align}
\end{definition}

The sign in \eqref{gdu-du} determines the parameter correspondence below. For arbitrary $r,s$, the monomials $u^i h^j d^k$ with $i,j,k\in\mathbb{N}$ form a $\Bbbk$-basis, and $\operatorname{GKdim}L=3$ \cite[Theorem 2.1 and Corollary 2.4]{CassidyShelton2004}. In particular, $\Bbbk[h]$ embeds in $L$ even if $rs=0$.

\begin{proposition}[{\cite[Lemma 1.1]{CarvalhoLopes2009}}]\label{parameter-correspondence}
Let $r,s,\varepsilon,\mu\in\Bbbk$ and $\lambda\in\Bbbk^*$. There is an isomorphism
\begin{equation}\label{linear-correspondence}
L(\lambda x + \mu,r,s,\varepsilon) \cong A(r+s,-rs,\lambda\varepsilon+(r-1)\mu)
\end{equation}
that preserves the generators $d,u$ and sends
\[
h\longmapsto\lambda^{-1}(sud-du-\mu).
\]
In particular, if $x^2-\alpha x -\beta$ splits over $\Bbbk$ with roots $r,s\in\Bbbk$, then
\begin{equation}\label{du-gdu-correspondence}
A(\alpha,\beta,\gamma)\cong L(x, r, s, \gamma), \  \alpha = r + s, \quad \beta = -rs.
\end{equation}
\end{proposition}

This is the corrected correspondence recorded in \cite[Lemma 1.1]{CarvalhoLopes2009}. It applies even if one of $r,s$ is zero. For example, eliminating $h$ from \eqref{gdu-dh}--\eqref{gdu-du} when $f(X)=X$ gives
\[
d^2u=(r+s)dud-rs\,ud^2+\varepsilon d,\  du^2=(r+s)udu-rs\,u^2d+\varepsilon u,
\]
which agrees with \eqref{du-first}--\eqref{du-second}.

\begin{proposition}[{\cite[Propositions 2.5 and 2.6]{CassidyShelton2004}; \cite[Proposition 1.2]{CarvalhoLopes2009}}]\label{gdu-noetherian}
For $L = L(f,r,s,\varepsilon)$, the following are equivalent: 
\begin{enumerate}[label=\textup{(\arabic*)},leftmargin=*]
\item $L$ is Noetherian;
\item $L$ is a domain;
\item $rs\ne0$.
\end{enumerate}
\end{proposition}

In the case $rs \neq 0$, the algebra has an explicit iterated Ore presentation \cite{Ore1933}. Recall that $R[z;\sigma,\delta]$ is free as a left $R$-module on the set $\{z^j\mid j\in\mathbb{N}\}$, and its multiplication is determined by $
za = \sigma(a)z+\delta(a)$ for all $a\in R$, where $\delta$ is a {\em left} $\sigma$-{\em derivation}, that is, $\delta(ab)=\sigma(a)\delta(b)+\delta(a)b$. 

\begin{proposition}[{\cite[Lemma 1.4]{CarvalhoLopes2009}}]\label{gdu-ore}
Assume $rs\ne0$. Let $\sigma$ be the automorphism of $\Bbbk[h]$ defined by $\sigma(h)=rh-\varepsilon$. Then
\begin{equation}\label{gdu-ore1}
L(f,r,s,\varepsilon)\cong\Bbbk[h][d;\sigma][u;\tau,\delta],
\end{equation}
where the automorphism $\tau$ and the $\tau$-derivation $\delta$ of $\Bbbk[h][d;\sigma]$ are determined by
\begin{align*}
\tau(h)&=r^{-1}(h+\varepsilon),&\  \tau(d)&=s^{-1}d,\\
\delta(h)&=0,&\delta(d)&=s^{-1}f(h).
\end{align*}
Consequently, $
\left\{h^i d^j u^k\mid i,j,k\in\mathbb{N}\right\}$ is a $\Bbbk$-basis of $L(f,r,s,\varepsilon)$.
\end{proposition}

The Ore presentation also gives a direct growth computation. Choose an integer $N\geq1$ with $2N>\deg f$ when $f\ne0$, and assign weights $1,N,N$ to $h,d,u$, respectively. Let $F_m$ be the span of the PBW monomials $h^i d^j u^k$ of weight $i+N(j+k)\leq m$. Reordering a word into PBW form does not increase its weight, so $(F_m)_{m\geq0}$ is the filtration by weighted word length. The terms involving $\varepsilon$ and $f$ have strictly lower weight than the corresponding quadratic terms. Thus the associated graded relations are
\[
\bar d\bar h=r\bar h\bar d,\quad
\bar h\bar u=r\bar u\bar h,\quad {\rm and} \quad \bar d\bar u=s\bar u\bar d.
\]
The ordinary total degree of a PBW monomial is at most its weight and at least its weight divided by $N$. Counting triples of total degree at most $\lfloor m/N\rfloor$ and at most $m$ therefore gives
\[
\binom{\lfloor m/N\rfloor+3}{3}\leq \dim_{\Bbbk}(F_m)\leq\binom{m+3}{3}.
\]
For $V=\operatorname{span}_{\Bbbk}\{1,h,d,u\}$ one has $F_m\subseteq V^m\subseteq F_{Nm}$. These inclusions give cubic growth of $\dim_{\Bbbk}V^m$, so $\operatorname{GKdim}L=3$.

There is a second presentation that is often useful when studying automorphisms and skew derivations. For a commutative algebra $B$, an automorphism $\rho$ of $B$ and an element $t\in B$, the generalized Weyl algebra $B(\rho,t)$ has generators $B,x,y$ with relations
\[
xg=\rho(g)x,\qquad gy=y\rho(g)\quad {\rm for} \ g\in B,\qquad
yx=t,\qquad xy=\rho(t).
\]
For $L$ with $rs\ne0$, take $B=\Bbbk[h,t]$, a polynomial algebra in two variables, and set
\begin{equation}\label{gdu-gwa}
\rho(h)=rh-\varepsilon,\qquad \rho(t)=st-f(h).
\end{equation}
Then
\[
L(f,r,s,\varepsilon)\cong B(\rho,t),\quad {\rm with} \quad d\longmapsto x,\quad u\longmapsto y,
\]
where $t$ corresponds to $ud$ \cite[Lemma 2.7]{CassidyShelton2004} and \cite[Lemma 1.3]{CarvalhoLopes2009}.

\begin{definition}[{\cite[Section 1.3]{CarvalhoLopes2009}}]\label{conformal}
A generalized down-up algebra $L(f, r, s, \varepsilon)$ is {\em conformal} if there exists $g\in\Bbbk[x]$ such that
\[
f(x) = sg(x) - g(rx - \varepsilon).
\]
\end{definition}

In the Noetherian case, a polynomial $g$ as above gives the normal element
\[
z=du-g(rh-\varepsilon)=s\bigl(ud-g(h)\bigr),
\]
which satisfies $zh=hz$, $dz=szd$ and $zu=suz$ \cite[Section 1.3]{CarvalhoLopes2009}. Note that conformality is an additional property of the polynomial data; it is not imposed in Definition \ref{generalized-down-up1}.

\begin{example}\label{gdu-examples}
\begin{enumerate}
\item [\rm (i)] \cite[Introduction and Section 3.1]{LaunoisLopes2013} \label{gdu-quantum-affine} 
The specialization $f=0$ and $\varepsilon=0$ in Definition \ref{generalized-down-up1}, with $r,s\in\Bbbk^*$, gives
\[
L(0,r,s,0) = \Bbbk\{ d,u,h\} / \langle dh-rhd,\ hu-ruh,\ du-sud\rangle.
\]
This specialization of the Cassidy and Shelton family is a three variable quantum affine space. In particular, $L(0,1,1,0)$ is the commutative polynomial algebra $\Bbbk[h, d, u]$, and Example \ref{polynomial-calculus} supplies its connected integrable three-dimensional calculus. The vanishing polynomial also makes this family conformal, by taking $g=0$.

\item [\rm (ii)] \cite[Introduction]{Smith1990}; \cite[Section 7]{LaunoisLopes2013} \label{Smith}
For a polynomial $P\in\Bbbk[x]$, consider the presentation
\[
[H,E]=E,\  [H,F]=-F \quad {\rm and} \quad [E,F]=P(H).
\]
These are the algebras similar to $U(\mathfrak{sl}_2)$ studied by Smith \cite{Smith1990}. Under the identification $h=H$, $d=E$ and $u=F$, they have the generalized down-up presentation $ L(-P,1,1,1)$. Indeed, \eqref{gdu-dh} and \eqref{gdu-hu} give the first two commutators, while \eqref{gdu-du} gives the third. For instance, $P(x) = x^2$ yields the nonlinear example $L(-x^2,1,1,1)$. The polynomial has degree two, so Proposition \ref{parameter-correspondence} does not eliminate $h$ by its affine linear substitution.

\item [\rm (iii)] \cite[Lemma 1.1]{CarvalhoLopes2009} \label{gdu-enveloping} The correspondence with ordinary down-up algebras includes the two enveloping algebras discussed above:
\[
L(x,1,1,0)\cong U(\mathfrak h_3),\ L(x,1,1,-2)\cong U(\mathfrak{sl}_2(\Bbbk)).
\]
In both cases the distinguished generator $h$ corresponds to $ud-du$. In the second isomorphism, with $d=e$ and $u=f$ as in Example \ref{sl-two}, it therefore corresponds to $-H$.
\end{enumerate}
\end{example}

\section{Differential and integral calculi on down-up algebras}\label{Differentialandintegralcalculusdownup}

In this section, we construct three-dimensional differential calculi on generalized down-up algebras and apply them to ordinary down-up algebras. As in the proof of \cite[Theorem 3.2]{RubianoReyes2026BQAThree}, we begin by specifying how the algebra acts on the one-forms and then check the Leibniz rule on the defining relations. We use $D$ for the exterior differential throughout this section in order to distinguish it from the down generator $d$. Products of forms are written by juxtaposition.

Fix some notation: for a nonzero scalar $a\in \Bbbk$ and $n\geq1$, let $[n]_a := 1+a+\cdots+a^{n-1}$, with $[0]_a := 0$. Over a field of characteristic zero, all the $[n]_a$, $n\geq1$, are nonzero precisely when $a=1$ or $a$ has infinite multiplicative order.

The following observation will be used to check integrability in both constructions.

\begin{lemma}\label{complementary-forms}
Suppose a differential calculus $(\Omega S, d)$ over an algebra $S$ has left and right bases in every degree $k$ {\rm (}$0 \le k \le n${\rm )} consisting of the increasing products $e_I = e_{i_1} \wedge \cdots \wedge e_{i_k}$, for $I = \{i_1 < \cdots < i_k\} \subseteq \{1,\dots,n\}$ {\rm (}with $e_\emptyset = 1${\rm )}, of a fixed list of one-forms $e_1,\dots,e_n$. Assume $e_i \wedge e_i = 0$ and $e_j \wedge e_i = -q_{ij} e_i \wedge e_j$ for $i < j$, with $q_{ij} \in \Bbbk^*$. Then $\omega = e_1 \wedge \cdots \wedge e_n$ is an integrating volume form.
\end{lemma}
\begin{proof}
The form $\omega$ is a free generator of $\Omega^n S$ on both sides. If $|I|=k$ and $|J|=n-k$, then the relations $e_i \wedge e_i = 0$ and $e_j \wedge e_i = -q_{ij} e_i \wedge e_j$ ($i < j$, $q_{ij} \in \Bbbk^*$) imply that $e_I \wedge e_J = 0$ when $J \ne I^c$, and $e_I \wedge e_{I^c} = c_I \omega$ for some $c_I \in \Bbbk^*$. Write $\theta = \sum_{|I|=k} a_I e_I$ using the left basis. Evaluation on the right basis gives
\[
\ell_{\pi_\omega}^{k}(\theta)(e_J) = \pi_\omega(\theta \wedge e_J) = c_{J^c}\nu_\omega(a_{J^c}).
\]
The values on the $e_J$ determine a right $S$-linear functional, and every list of such values determines a unique list of the $a_I$. Thus $\ell_{\pi_\omega}^{k}$ is bijective for every $1 \le k \le n-1$, and Proposition~\ref{integrability} now gives the result.
\end{proof}

\subsection{Diagonal calculi for generalized down-up algebras}\label{gdu-diagonal}

Let $L=L(f,r,s,\varepsilon)$ and write $T(h)=rh-\varepsilon$. If $r\ne0$, the inverse substitution is $T^{-1}(h)=r^{-1}(h+\varepsilon)$. Write $D$ for the differential and put $\xi=D(u)$, $\eta=D(h)$, and $\zeta=D(d)$. A diagonal left action on the free right $L$-module with basis $\xi,\eta,\zeta$ is one of the form
\[
a\xi=\xi\nu_u(a),\quad a\eta=\eta\nu_h(a),\quad a\zeta=\zeta\nu_d(a),\quad a\in L,
\]
where $\nu_u,\nu_h,\nu_d$ are algebra automorphisms of $L$.

\begin{proposition}\label{gdu-diagonal-obstruction}
If a generalized down-up algebra $L$ admits a connected differential calculus with a free right basis $D(u),D(h),D(d)$ for its one-forms and with diagonal left action, then $rs\ne0$ and $f$ is constant. For $rs\ne0$, the assignments
\begin{align}
\nu_u(u)&=u,&\nu_u(h)&=T(h),&\nu_u(d)&=sd,\notag\\
\nu_h(u)&=r^{-1}u,&\nu_h(h)&=h,&\nu_h(d)&=rd,\label{gdu-diagonal-fixed}\\
\nu_d(u)&=s^{-1}u,&\nu_d(h)&=T^{-1}(h),&\nu_d(d)&=d\notag
\end{align}
define automorphisms of $L$ if and only if $f(T(h))=sf(h)$. In characteristic zero, a connected three-dimensional integrable calculus with the assignments exactly as in \eqref{gdu-diagonal-fixed} -- in particular with $\nu_u(u)=u$ and $\nu_d(d)=d$ --, exists precisely when $f=c\in\Bbbk$ and $(s-1)c=0$.
\end{proposition}
\begin{proof}
Since $h\eta=\eta\nu_h(h)$, the diagonal action and the Leibniz rule give $D(f(h))\in\eta L$. The coefficient of $\eta$ in $D(hu-ruh+\varepsilon u)$ is $u-r\nu_h(u)$. If $r=0$, this contradicts the PBW basis, which contains $u$.
Differentiating \eqref{gdu-du} gives
\[
\xi\bigl(\nu_u(d)-sd\bigr)+\zeta\bigl(u-s\nu_d(u)\bigr)+D(f(h))=0.
\]
Freeness gives $D(f(h))=0$ and $u=s\nu_d(u)$. The latter identity forces $s\ne0$. If the calculus is connected, then
$f(h)\in\Bbbk$. By the PBW basis of $L$, $h$ generates a polynomial subalgebra, so $f$ is constant.

Each assignment in \eqref{gdu-diagonal-fixed} preserves \eqref{gdu-dh} and \eqref{gdu-hu}. On the relation $du-sud+f(h)=0$, the images under $\nu_u,\nu_h,\nu_d$ are, respectively, $f(T(h))-sf(h)$, zero, and $f(T^{-1}(h))-s^{-1}f(h)$. The two polynomial conditions are equivalent, and the inverse assignments preserve the relations under the same condition. Thus the assignments define automorphisms exactly when $f(T(h))=sf(h)$. Since $\nu_h(h)=h$, the remaining first-order compatibility condition is $f'(h)=0$. In characteristic zero this is equivalent to $f=c\in\Bbbk$. The automorphism condition then becomes $(s-1)c=0$. For these parameters, the construction in Theorem \ref{gdu-constant} supplies the required extension with both own-generator scalings equal to $1$. For a nonzero constant and $s\ne1$, that theorem uses rescaled assignments under its quantum-integer hypothesis.
\end{proof}

\begin{theorem}\label{gdu-constant}
Assume that $\operatorname{char}\Bbbk=0$ and that $r,s\in\Bbbk^*$. Let $c\in\Bbbk$. If $c=0$, then $L(c,r,s,\varepsilon)$ is differentially smooth. If $c\ne0$, the same conclusion holds whenever $[n]_s\ne0$ for every $n\geq1$. In characteristic zero the latter condition means that $s=1$ or $s$ has infinite multiplicative order.
\end{theorem}
\begin{proof}
Set $L=L(c,r,s,\varepsilon)$ and let $a_u=a_d=1$ if $c=0$. If $c\ne0$, let $a_u=s^{-1}$ and $a_d=s$. Keep the other assignments in \eqref{gdu-diagonal-fixed} and set $\nu_u(u)=a_u u$ and $\nu_d(d)=a_d d$. When $c\ne0$ and $s\ne1$, these two assignments differ from those in Proposition \ref{gdu-diagonal-obstruction}. The relations \eqref{gdu-dh} and \eqref{gdu-hu} remain invariant up to invertible scalar factors. The images of $du-sud+c$ under $\nu_u,\nu_h,\nu_d$ are, respectively, $(1-sa_u)c$, zero, and $(1-s^{-1}a_d)c$ in $L$. Our choices make these expressions zero. The inverse assignments also preserve the relations, and the resulting automorphisms commute. Indeed, on $h$ they act by $T$, the identity, and $T^{-1}$, and on $u,d$ they act by scalars.

We use the diagonal left action above and impose the exterior relations
\begin{equation}\label{gdu-constant-exterior}
\xi^2=\eta^2=\zeta^2=0,\quad\eta\xi=-r\xi\eta,\qquad \zeta\xi=-s\xi\zeta,\quad\zeta\eta=-r\eta\zeta.
\end{equation}
To see that the ordered products of distinct forms give a free right $L$-basis, let $E$ be the graded scalar algebra given by \eqref{gdu-constant-exterior}. Its ordered square-free products $e_I$, indexed by subsets $I\subseteq\{1,2,3\}$ in the order $\xi,\eta,\zeta$, are a basis. The scalar swaps give the same coefficient in either order at each overlap, and a repeated form reduces to zero. On $E\otimes_{\Bbbk}L$ define $(e_I a)(e_J b)=(e_Ie_J)\nu_J(a)b$, where $\nu_J$ is the composition of the automorphisms corresponding to the factors of $e_J$ (which is independent of order since the automorphisms commute) and $\nu_{\varnothing}$ is the identity. This multiplication is associative because the automorphisms commute and $\nu_{J\cup K}=\nu_K\circ\nu_J$ when $J\cap K=\varnothing$. If some form is repeated in a triple product, both associations vanish. The multiplication realizes the diagonal commutations and \eqref{gdu-constant-exterior}. Moving the algebra generators to the right then proves that these relations present $E\otimes_{\Bbbk}L$.

Define $D(u)=\xi$, $D(h)=\eta$, $D(d)=\zeta$ and $D(\xi)=D(\eta)=D(\zeta)=0$. Applying the graded Leibniz rule to
the three algebra relations gives
\begin{align*}
D(hu-uT(h))&=\eta u+\xi T(h)-\xi T(h)-r(r^{-1}\eta u)=0,\\
D(dh-T(h)d)&=\zeta h+r\eta d-r\eta d-\zeta T(T^{-1}(h))=0,\\
D(du-sud+c)&=\zeta u+s\xi d-s\xi d-s(s^{-1}\zeta u)=0.
\end{align*}
The differential of each mixed commutation vanishes by \eqref{gdu-constant-exterior}. For instance, $D(h\xi-\xi T(h))=\eta\xi+r\xi\eta$ and $D(d\xi-s\xi d)=\zeta\xi+s\xi\zeta$. The own-generator commutations give multiples of squares. The exterior relations have zero differential. This defines a differential graded algebra with $D^2=0$, and its ordered basis forms are products of differentials of algebra generators.

By Proposition \ref{gdu-ore}, the monomials $h^i d^j u^k$ form a basis of $L$. Write $D(v)=\xi\partial_u(v)+
\eta\partial_h(v)+\zeta\partial_d(v)$ for $v\in L$. For the basis monomials the coefficients are
\begin{align*}
\partial_h(h^i d^j u^k)&=i h^{i-1}d^j u^k,\\
\partial_d(h^i d^j u^k)&=[j]_{a_d} T^{-1}(h)^i d^{j-1}u^k,\\
\partial_u(h^i d^j u^k)&=s^j[k]_{a_u} T(h)^i d^j u^{k-1}.
\end{align*}
Terms with the corresponding exponent zero are interpreted as zero. If $D(v)=0$, the first equation and characteristic zero give $v\in\operatorname{span}_{\Bbbk}\{d^j u^k\}$. The second gives $v\in\Bbbk[u]$ and the third gives $v\in\Bbbk$, since all positive $[j]_{a_d}$ and $[k]_{a_u}$ are nonzero under our hypotheses. For $c\ne0$ we use $[n]_{s^{-1}}=s^{1-n}[n]_s$.

The ordered products of distinct forms are also a free left basis, since each $\nu_J$ is invertible. The top form
$\omega=\xi\eta\zeta$ is a volume form with $\nu_\omega=\nu_d\circ\nu_h\circ\nu_u$. Every pair of complementary basis forms has product a nonzero scalar multiple of $\omega$, and products of noncomplementary basis forms of
total degree three vanish. Lemma \ref{complementary-forms} therefore gives integrability. Since $\operatorname{GKdim} L=3$, the connected three-dimensional calculus proves differential smoothness.
\end{proof}

\begin{remark}\label{rem:gdu-constant-root-of-unity}
Suppose $c\ne0$ and $s$ has finite order $m>1$. For a diagonal calculus with $\nu_u(u)=a_u u+b_u$ ($a_u,b_u\in\Bbbk$), differentiating \eqref{gdu-du} forces $\nu_u(d)=sd$. Applying $\nu_u$ to $du-sud+c=0$ then gives
\[
(1-sa_u)c+s(1-s)b_u d=0.
\]
The PBW basis implies $a_u=s^{-1}$ and $b_u=0$. Hence $D(u^m)=[m]_{s^{-1}}\xi u^{m-1}=0$, although $u^m$ is a
nonscalar PBW monomial. In particular, the quantum integer condition is necessary for connectedness of the scaled calculus in Theorem \ref{gdu-constant}.
\end{remark}

\subsection{A calculus with triangular commutation rules}\label{sec:triangular-gdu}

In what follows, affine substitutions such as $T(h)$ and $G(h)$ are also used as automorphisms of $\Bbbk[h]$.

\begin{theorem}\label{triangular-gdu}
Let $\Bbbk$ have characteristic zero and let $L=L(f,r,s,\varepsilon)$ with $r,s\in\Bbbk^*$. Write $T(h)=rh-\varepsilon$. Suppose that there exist $\lambda\in\Bbbk^*$ and an affine automorphism $G(h)=qh+b$ of $\Bbbk[h]$ such that
\begin{equation}\label{eigen-gdu}
G\circ T=T\circ G,\quad f(T(h))=s\lambda f(h),\quad f(G(h))=\lambda f(h).
\end{equation}
If each of $\lambda$ and $q$ is either $1$ or has infinite multiplicative order, then $L$ is differentially smooth.
\end{theorem}

For $T(h)=rh-\varepsilon$ and $G(h)=qh+b$, the first condition in \eqref{eigen-gdu} is equivalent to $b(1-r)=\varepsilon(q-1)$, by comparing the constant terms of the two affine compositions.

\begin{proof}
Let $\kappa=\lambda/r$ and $\chi=\lambda-1$. For $g\in\Bbbk[h]$, set
\[
\partial_G g(h)=
\begin{cases}
\dfrac{g(G(h))-g(h)}{G(h)-h},& G\ne\mathrm{id},\\[5pt]
g'(h),&G=\mathrm{id}.
\end{cases}
\]
The divided difference is a polynomial. It satisfies $\partial_G(g_1g_2)=g_1(G(h))\partial_Gg_2+(\partial_Gg_1)g_2$.
Writing $F=\partial_G f$, we obtain from \eqref{eigen-gdu}
\begin{equation}\label{gdu-divided-identities}
(G(h)-h)F(h)=\chi f(h),\qquad F(T(h))=s\kappa F(h).
\end{equation}
For the second identity, the numerator of the divided difference at $T(h)$ is $s\lambda$ times its value at $h$, and the denominator is multiplied by $r$. When $G=\mathrm{id}$, the same identity follows by differentiating $f(T(h))=s\lambda f(h)$.

We use the order $u,h,d$ on the algebra generators. Set $R=\Bbbk[h][u;T^{-1}]$, so $hu=uT(h)$. On $R$ consider the right module freely generated by $\xi$ and $\eta$, with left action given by
\begin{equation}\label{gdu-base-commutations}
\begin{aligned}
u\xi&=\lambda\xi u, & h\xi&=\xi T(h)+\chi\eta u,\\
u\eta&=\kappa\eta u, & h\eta&=\eta G(h).
\end{aligned}
\end{equation}
To check the left action, apply $hu-uT(h)$ to both basis elements. On $\xi$ the difference in the $\xi$ coefficient is
$\lambda(T(h)u-uT(T(h)))=0$, and the $\eta$ coefficient is $\chi(\lambda-r\kappa)u^2=0$. On $\eta$ the difference vanishes by $G\circ T=T\circ G$ and $g(h)u=u g(T(h))$ for $g\in\Bbbk[h]$. Thus \eqref{gdu-base-commutations} defines an $R$-bimodule.

Write $M=\xi R\oplus\eta R$ for this bimodule and take its tensor algebra $T_R(M)$. We impose the exterior relations
\begin{equation}\label{gdu-base-exterior}
\xi^2=\eta^2=0,\qquad \eta\xi=-\kappa^{-1}\xi\eta.
\end{equation}
Let $W$ be the right span in $M\otimes_R M$ of $\xi^2$, $\eta^2$ and $w=\eta\xi+\kappa^{-1}\xi\eta$. Left multiplication by $u$ preserves $W$, since it multiplies the three generators by $\lambda^2$, $\kappa^2$ and $\lambda\kappa$, respectively, with a factor $u$ on the right. For $h$, the commutations in \eqref{gdu-base-commutations} give
\begin{align*}
h\xi^2&=\xi^2 T(T(h))+\chi\lambda w u,\\
h\eta^2&=\eta^2 G(G(h)),\\
hw&=w G(T(h))+\chi(q+1)\eta^2u.
\end{align*}
Here $T\circ G=G\circ T$ is used in the third identity. Hence $W$ is an $R$-subbimodule.

The length-$n$ words in $\xi,\eta$ are a right $R$-basis of $M^{\otimes_R n}$, by induction from $M=(\Bbbk\xi\oplus\Bbbk\eta)\otimes_{\Bbbk}R$. For $n\geq2$, the degree-$n$ part of the ideal generated by $W$ is the sum of $M^{\otimes_R i} W M^{\otimes_R j}$ over $i+j=n-2$. This sum is the right span of words with one of the three quadratic relators inserted. Indeed, coefficients to the left of the relator can be moved through it within $W$ and then through the remaining word. Conversely, each insertion lies in the ideal. Thus the quotient in degree $n$ is the corresponding degree of the scalar algebra with relations \eqref{gdu-base-exterior}, tensored over $\Bbbk$ with $R$. In degree two the three independent relators leave $\xi\eta$ as a basis. In higher degrees, each word has a repeated letter and reduces to zero by the scalar exchange and the square relations. Thus $\Omega R$ has right $R$-basis $1,\xi,\eta,\xi\eta$.

Define $D_R(u)=\xi$, $D_R(h)=\eta$, and $D_R(\xi)=D_R(\eta)=0$. The derivative of $hu-uT(h)$ is zero by \eqref{gdu-base-commutations}. Differentiating the mixed relations gives scalar multiples of the squares and of $r\xi\eta+\lambda\eta\xi$, all of which vanish by \eqref{gdu-base-exterior}. This defines $D_R$ on $\Omega R$ and $D_R^2=0$. In particular, $D_R(g(h))=\eta\partial_Gg(h)$.

We next extend the Ore presentation from degree zero to the forms. On $\Omega R$ define a degree-preserving automorphism $\Sigma$ by
\begin{equation}\label{gdu-sigma}
\Sigma(h)=T(h),\quad \Sigma(u)=su,\quad\Sigma(\eta)=r\eta,\quad \Sigma(\xi)=s\xi.
\end{equation}
Indeed, the substitution scales each defining relation, using $GT=TG$ in the relation for $h\eta$. The generators can be recovered from their images, and the inverse substitution with $T^{-1}$ and the reciprocal scalars also preserves the relations. Define a degree-zero $\Sigma$-derivation $\Delta$ by
\begin{equation}\label{gdu-delta}
\Delta(h)=\Delta(\eta)=0,\quad \Delta(u)=-f(h),\quad \Delta(\xi)=-\eta F(h).
\end{equation}
For its verification, induction in the degree of $g$ gives
\begin{equation}\label{gdu-commute-polynomial}
g(h)\xi=\xi g(T(h))+\chi\eta(\partial_G g(h))u \quad(g\in\Bbbk[h]).
\end{equation}
In the induction step, $uT(h)=hu$ converts the mixed term to the divided difference of $hg(h)$.

The twisted Leibniz rule gives $\Delta(hu-uT(h))=0$, since $T(h)$ commutes with $f(h)$. The relations $h\eta-\eta G(h)$ and $\eta^2$ are killed directly. For the other relations, the following computations use \eqref{gdu-commute-polynomial} and $u g(h)=g(T^{-1}(h))u$:
\begin{align*}
\Delta(h\xi-\xi T(h)-\chi\eta u) &=r\eta\bigl(\chi f(h)-(G(h)-h)F(h)\bigr),\\
\Delta(u\eta-\kappa\eta u) &=\eta\bigl(\lambda f(h)-f(G(h))\bigr),\\
\Delta(u\xi-\lambda\xi u) &=\xi\bigl(s\lambda f(h)-f(T(h))\bigr)
    +\eta\bigl(F(h)-s\kappa F(T^{-1}(h))\bigr)u,\\
\Delta(\xi^2) &=\xi\eta\bigl(\kappa^{-1}F(T(h))-sF(h)\bigr).
\end{align*}
Each expression vanishes by \eqref{eigen-gdu} and \eqref{gdu-divided-identities}. The image of $\eta\xi+\kappa^{-1}\xi\eta$ vanishes because both terms contain a repeated $\eta$ after moving the polynomial coefficient to the right. This proves that $\Delta$ is well defined.

On the generators, $D_R\Sigma=\Sigma D_R$ and $D_R\Delta=\Delta D_R$, and the respective product rules extend these identities to $\Omega R$. We also have
\begin{equation}\label{gdu-sigdelta}
\Sigma\Delta=\lambda\Delta\Sigma.
\end{equation}
For $u$ and $\xi$ this follows respectively from $f(T(h))=s\lambda f(h)$ and $rF(T(h))=s\lambda F(h)$. For $h$ and $\eta$ both sides vanish. The twisted product rules again extend the identity.

Let $C=\Omega R[d;\Sigma,\Delta]$. Its degree-zero subalgebra is $L$, since the Ore relations in degree zero are $dh=T(h)d$ and $du=sud-f(h)$. The relations also show that $u^ih^jd^k$, with $i,j,k\in\mathbb{N}$, form a basis. For homogeneous $\theta\in\Omega R$, the assignments
\[
\tau(\theta)=(-1)^{|\theta|}\Sigma(\theta),\qquad \tau(d)=\lambda d
\]
give a graded automorphism of $C$. In fact, applying $\tau$ to the Ore relation $d\theta=\Sigma(\theta)d+\Delta(\theta)$ gives that relation at $\Sigma(\theta)$, using \eqref{gdu-sigdelta}. The inverse uses the reciprocal scalar and the inverse of the signed automorphism of $\Omega R$. Set
\begin{equation}\label{gdu-odd-extension}
\Omega L=C[\zeta;\tau]/\langle\zeta^2\rangle,\qquad |\zeta| = 1.
\end{equation}
The Ore normal form in $C[\zeta;\tau]$ is a finite sum $\sum_j c_j\zeta^j$ with $c_j\in C$. Since $\zeta^2c=\tau^2(c)\zeta^2$ for $c\in C$, the span $\bigoplus_{j\geq2}C\zeta^j$ is a two-sided ideal. Each of its terms is a left multiple of $\zeta^2$, so it is exactly the ideal in \eqref{gdu-odd-extension}. The quotient is therefore $C\oplus C\zeta$ as a left $C$-module. The new relations in $\Omega L$ are
\begin{equation}\label{gdu-triangular-commutations}
\begin{aligned}
d\xi&=s\xi d-\eta F(h), & d\eta&=r\eta d, & d\zeta&=\lambda^{-1}\zeta d,\\
u\zeta&=s^{-1}\zeta u, & h\zeta&=\zeta T^{-1}(h), & \zeta^2&=0,\\
\zeta\xi&=-s\xi\zeta, & \zeta\eta&=-r\eta\zeta.
\end{aligned}
\end{equation}
The first formula gives the term needed when differentiating $du-sud+f(h)$. This term also explains the use of $h\xi=\xi T(h)+\chi\eta u$ in the calculus on $R$.

Extend $D_R$ by $D(d)=\zeta$ and $D(\zeta)=0$. Differentiating the Ore relation in $C$ gives zero because $D_R$ commutes with $\Sigma,\Delta$ and $\zeta\theta=(-1)^{|\theta|}\Sigma(\theta)\zeta$. The latter relation is stable under $D$ because $D_R$ raises degree by one. The derivative of $\zeta d-\lambda d\zeta$ is a multiple of $\zeta^2$. Thus $D$ is well defined and $D^2=0$.

The Ore normal form in $C=\Omega R[d;\Sigma,\Delta]$ is a finite sum $\sum_j\theta_jd^j$ with $\theta_j\in\Omega R$. Expanding each $\theta_j$ in the right $R$-basis above gives
\[
C=\bigoplus_{w\in\{1,\xi,\eta,\xi\eta\}}wL.
\]
The sum is direct by uniqueness of the two expansions. Since $\tau$ restricts to an automorphism of $L$, the rule $a\zeta=\zeta\tau^{-1}(a)$ for $a\in L$ identifies $wL\zeta$ with $w\zeta L$, with unique right coefficients. Together with $\Omega L=C\oplus C\zeta$, this gives a right $L$-basis of increasing products of $\xi,\eta,\zeta$.

To obtain the left bases, consider the maps
\begin{equation}\label{gdu-diagonal-maps}
\begin{array}{c|ccc}
 &u&h&d\\ \hline
\nu_\xi&\lambda u&T(h)&s d\\
\nu_\eta&\kappa u&G(h)&r d\\
\nu_\zeta&s^{-1}u&T^{-1}(h)&\lambda^{-1}d
\end{array}
\end{equation}
These are automorphisms of $L$ by \eqref{eigen-gdu}, with inverses given by the inverse affine substitutions and reciprocal scalars. The commutation rules for the algebra generators extend by multiplication to
\[
a\xi=\xi\nu_\xi(a)+\eta\vartheta(a),\qquad
a\eta=\eta\nu_\eta(a),\qquad
a\zeta=\zeta\nu_\zeta(a)
\quad(a\in L),
\]
where $\vartheta(a)$ is the right coefficient of $\eta$ in $a\xi$. For example,
\[
a\xi\zeta
=\xi\zeta(\nu_\zeta\circ\nu_\xi)(a)+\eta\zeta\nu_\zeta(\vartheta(a)).
\]
In $a\xi\eta$ and $a\xi\eta\zeta$, the term containing $\vartheta(a)$ vanishes by $\eta^2=0$. Thus in the ordered bases $(\xi,\eta,\zeta)$ and $(\xi\eta,\xi\zeta,\eta\zeta)$ the only off-diagonal contributions go from $\xi$ to $\eta$ and from $\xi\zeta$ to $\eta\zeta$, respectively. In every degree, the diagonal map for a word is the composition of the automorphisms for its successive factors. Each diagonal map is invertible. Given the right coefficients, one can therefore solve successively for the left coefficients, starting with the first basis word in each order. This gives a unique expansion on the left, and the left bases consist of products of exact forms. Hence the density condition holds.

To verify connectedness, the Leibniz rule on an ordered monomial gives
\begin{align*}
D(u^ih^jd^k) ={}&\xi[i]_\lambda u^{i-1}h^jd^k +\eta\,\kappa^i u^i\partial_G(h^j)d^k\\
 &+\zeta\,s^{-i}[k]_{\lambda^{-1}}u^i\bigl(T^{-1}(h)\bigr)^j d^{k-1},
\end{align*}
where a term is omitted when its lowered exponent is negative. The $\zeta$ coefficient of the differential of a linear combination can vanish only if all terms of positive $d$-degree vanish. Indeed, $T^{-1}$ is invertible and $[k]_{\lambda^{-1}}=\lambda^{1-k}[k]_\lambda$ for $k\geq1$. The leading coefficient of $\partial_G(h^j)$ is $[j]_q$ in degree $j-1$ for $j\geq1$. The $\eta$ coefficient therefore eliminates the positive powers of $h$. Finally, the $\xi$ coefficient eliminates the positive powers of $u$ because $[i]_\lambda\ne0$ for $i\geq1$. Thus $\ker(D|_L)=\Bbbk 1$.

The form $\omega=\xi\eta\zeta$ freely generates the top forms on both sides. The terms that replace $\xi$ by $\eta$ vanish on $\omega$. In particular, its twisting automorphism is given by
\[
\nu_\omega(u)=\frac{\lambda^2}{rs}u,\quad \nu_\omega(h)=G(h),\quad \nu_\omega(d)=\frac{rs}{\lambda}d.
\]
The scalar exchanges in \eqref{gdu-base-exterior} and \eqref{gdu-triangular-commutations} now give integrability by Lemma \ref{complementary-forms}. Since $\operatorname{GKdim}L=3$ by Proposition \ref{gdu-ore} and its growth computation, this calculus establishes differential smoothness.
\end{proof}

\begin{corollary}\label{cor:gdu-eigen}
Let $r,s\in\Bbbk^*$ and let $\Bbbk$ have characteristic zero. If $f(rh-\varepsilon)=s f(h)$, then $L(f,r,s,\varepsilon)$ is differentially smooth.
\end{corollary}

\begin{proof}
For $G=\mathrm{id}$, the last condition in \eqref{eigen-gdu} reads $f=\lambda f$. It forces $\lambda=1$ when $f\ne0$, and we can also choose $\lambda=1$ when $f=0$. With this choice, the other polynomial condition is precisely $f(T(h))=s f(h)$. Theorem \ref{triangular-gdu} applies with $F=f'$ and with a single mixed term given by $d\xi=s\xi d-\eta f'(h)$.
\end{proof}

\subsection{Ordinary down-up algebras}

The correspondence in Proposition \ref{parameter-correspondence} gives the following application of Theorem \ref{triangular-gdu}. Throughout this subsection, the ground field has characteristic zero.

\begin{theorem}\label{TheoremDU}
Let $\beta\ne0$ and suppose that $x^2-\alpha x - \beta$ splits over $\Bbbk$, with roots $r,s$. If $r/s=1$ or $r/s$ has infinite multiplicative order, then $A(\alpha,\beta,0)$ is differentially smooth.
\end{theorem}
\begin{proof}
The roots satisfy $r+s=\alpha$ and $rs=-\beta$, so both are nonzero. Put $h=sud-du$. Proposition \ref{parameter-correspondence} identifies $A(\alpha,\beta,0)$ with $L(X,r,s,0)$, with $h$ as the distinguished generator. The affine map $T(h)=rh$ commutes with $G(h)=\lambda h$, where $\lambda=r/s$. For $f(x) = x$, we have
\[
f(T(h))=rh=s\lambda f(h),\qquad f(G(h))=\lambda f(h).
\]
Both $\lambda$ and the slope of $G$ satisfy the order hypothesis in Theorem \ref{triangular-gdu}. The connected integrable calculus furnished there therefore transfers along the isomorphism.
\end{proof}

The condition in Theorem \ref{TheoremDU} is independent of the order of the roots. It includes $A(2r,-r^2,0)$ for every $r\in\Bbbk^*$, and in particular $A(2,-1,0)\cong U(\mathfrak h_3)$. In the calculus of Theorem \ref{triangular-gdu}, the form $D(h)$ occurs in the rule
\[
dD(u)=sD(u)d-D(h).
\]
Thus $D(h)$ can belong to the bimodule generated by $D(d)$ and $D(u)$ while forming a right basis together with them. This distinction between bimodule generators and a right basis is important for the following comparison.

\subsection{The diagonal rank-two construction}

We examine the construction with independent right generators $D(d)$ and $D(u)$ and separate scaling automorphisms for the two forms. Suppose that $\operatorname{char}\Bbbk=0$, $\beta\ne0$, and a square root of $\alpha^2+4\beta$ lies in $\Bbbk$. Choose
\[
a_0=\frac{-\alpha+\sqrt{\alpha^2+4\beta}}{2\beta}, \quad b_0=\frac{\alpha\pm\sqrt{\alpha^2+4\beta}}{2}.
\]
The sign in the definition of $b_0$ can be chosen independently of the choice of square root. Put
\[
p_d=(\alpha a_0-1)^{-1},\qquad p_u=\frac{\alpha-b_0}{b_0}.
\]
These scalars are nonzero since
\[
\beta a_0^2+\alpha a_0-1=0,\qquad b_0^2-\alpha b_0-\beta=0.
\]
On $R=A(\alpha,\beta,\gamma)$, consider the assignments
\begin{align}
\nu_d(d)&=p_d d,& \nu_d(u)&=a_0 u,&
\nu_u(d)&=b_0 d,& \nu_u(u)&=p_u u.\label{du-diagonal-autos}
\end{align}
When they define automorphisms, set $\Omega^1 R=\xi R\oplus\eta R$, with $x\xi=\xi\nu_d(x)$ and $x\eta=\eta\nu_u(x)$ for $x\in R$, and prescribe $D(d)=\xi$ and $D(u)=\eta$.

\begin{proposition}\label{du-diagonal-obstruction}
Under the assumptions above with the prescriptions in \eqref{du-diagonal-autos} give a first-order calculus exactly when $\gamma=0$. In that case, let $r=a_0^{-1}$ and $s=\alpha-r$. This calculus is connected exactly when $b_0=s\ne r$ and $r/s$ has infinite multiplicative order. Every differential calculus extending this connected first-order calculus has zero forms in degrees at least two.
\end{proposition}

\begin{proof}
For $c_d,c_u\in\Bbbk^*$, the scaling $d\mapsto c_d d$, $u\mapsto c_u u$ preserves \eqref{du-first} and \eqref{du-second} exactly when $\gamma(c_dc_u-1)=0$. Indeed, substituting the scaled generators into the first relator yields $c_d\gamma(c_dc_u-1)d$, and the second yields $c_u\gamma(c_dc_u-1)u$. The inverse scaling then also preserves the relators. In particular, both scalings in \eqref{du-diagonal-autos} are automorphisms for $\gamma=0$.

Let $R_1,R_2$ be the left sides of \eqref{du-first} and \eqref{du-second} minus their right sides. For a prospective Leibniz map on the free algebra with values in the indicated $R$-bimodule, write $D(R_i)=\xi P_i+\eta Q_i$. The four coefficients are
\begin{align*}
P_1&=(1+p_d-\alpha p_d a_0)du-\bigl(\alpha+\beta a_0(1+p_d)\bigr)ud-\gamma,\\
Q_1&=(b_0^2-\alpha b_0-\beta)d^2,\\
P_2&=(1-\alpha a_0-\beta a_0^2)u^2,\\
Q_2&=\bigl(b_0(1+p_u)-\alpha\bigr)du
      -\bigl(\alpha b_0 p_u+\beta(1+p_u)\bigr)ud-\gamma.
\end{align*}
The two quadratic identities and the definitions of $p_d,p_u$ give $P_1=Q_2=-\gamma$ and $P_2=Q_1=0$. Consequently the map descends precisely when $\gamma=0$. When it descends, its image generates $\Omega^1 R$ because $D(d)=\xi$ and $D(u)=\eta$.

Now suppose $\gamma=0$. The roots of $x^2-\alpha x -\beta$ are $r=a_0^{-1}$ and $s=\alpha-r$, and $p_d=r/s$. If $b_0=r$, the element $h_{b_0}=b_0 ud-du$ satisfies
\[
D(h_{b_0})=(a_0 b_0-1)\xi u=0.
\]
By Proposition \ref{du-pbw}, $h_{b_0}$ is not a scalar. If $b_0=s\ne r$ and $p_d$ has finite multiplicative order $N$, then $p_d\ne1$, and
\[
D(d^N)=\xi[N]_{p_d} d^{N-1}=0.
\]
Again $d^N$ is not a scalar, so the calculus is disconnected in both cases.

It remains to consider $b_0=s\ne r$ and $p_d=r/s$ of infinite order. Here $p_u=p_d$. In the presentation with $h_{b_0}=b_0 ud-du$, the monomials $u^i h_{b_0}^j d^k$ form a PBW basis, and $h_{b_0}u=ru h_{b_0}$. Since $D(h_{b_0})$ and $D(d)$ belong to $\xi R$, the coefficient of $\eta$ in $D(u^i h_{b_0}^j d^k)$ for $i>0$ is
\[
[i]_{p_d} u^{i-1}h_{b_0}^j d^k.
\]
All $[i]_{p_d}$ are nonzero. Thus a vector in $\ker D$ is a linear combination of $h_{b_0}^j d^k$. Using $\nu_d(h_{b_0})=b_0^{-1}h_{b_0}$, we obtain
\[
D(h_{b_0}^j d^k)=\xi\Bigl((a_0 b_0-1)[j]_{p_d} u h_{b_0}^{j-1}d^k
                    + b_0^{-j}[k]_{p_d} h_{b_0}^j d^{k-1}\Bigr),
\]
where the first term is omitted if $j=0$ and the second if $k=0$. The two types of monomials have different $u$-exponents. Since $a_0 b_0\ne1$, PBW independence forces every coefficient with $j>0$ or $k>0$ to vanish in an element of $\ker D$. Hence $\ker D=\Bbbk$.

In any graded extension, applying $D$ to the commutation rules gives
\[
(1+p_d)\xi^2=(1+p_u)\eta^2=0,
\qquad \eta\xi=-a_0\xi\eta,\qquad \xi\eta=-b_0\eta\xi.
\]
Under connectedness, $p_d=p_u$ has infinite order and $a_0 b_0\ne1$, so all four products of two basis forms vanish. The diagonal left action moves every intermediate coefficient to the right. The density condition therefore gives $\Omega^2 R=0$ and $\Omega^n R=0$ for every $n\ge2$.
\end{proof}

In particular, the diagonal rank-two construction cannot furnish a connected calculus of the dimension required by Proposition \ref{du-pbw}. Theorem \ref{TheoremDU} uses the additional right basis form $D(h)$ and the mixed commutation rule displayed above. Note that the two quadratic identities governing the scalars in \eqref{du-diagonal-autos} are conditions for the first-order calculation. The quotient of the roots instead governs connectedness of the three-dimensional construction in Theorem \ref{TheoremDU}.

\section{Examples}\label{ExamplesDU}

We next illustrate the preceding results in the two families. All examples in this section are over a field of characteristic zero.

\subsection{Down-up algebras}

\begin{example}[The Heisenberg enveloping algebra]\label{heisenberg-smooth}
In Example \ref{Heisenberg}, $A(2,-1,0)\cong U(\mathfrak h_3)$. The polynomial $x^2 - 2x + 1$ has the double root $1$. In the notation of Theorem \ref{triangular-gdu}, $h=ud-du$, $f(x) = x$, $\lambda=1$, and $G=T=\mathrm{id}$. The forms anticommute and have zero squares. The sole nontrivial commutator of an algebra generator with a basis one-form is
\[
dD(u)-D(u)d=-D(h).
\]
For $i,j,k\in\mathbb N$, the differential of $u^ih^jd^k$ has the three ordinary exponent factors $i,j,k$. The form $D(u)D(h)D(d)$ is central and integrating. Heisenberg calculi with scalar commutators of generators and one-forms are described in \cite[Section 7]{MartiniPostKersten1995}. Theorem \ref{TheoremDU} also treats the double-root family $A(2r,-r^2,0)$ for every $r\in\Bbbk^*$.
\end{example}

\begin{example}[Two distinct roots]\label{distinct-down-up}
Consider $A(3,-2,0)$ and take $r=2$ and $s=1$. Then $h=ud-du$, $T(h)=G(h)=2h$, and $\lambda=2$. Theorem \ref{TheoremDU} applies. Here $hD(u)=2D(u)h+D(h)u$ and $dD(u)=D(u)d-D(h)$. Proposition \ref{du-diagonal-obstruction} gives a connected first-order calculus for the choice $a_0=1/2$, $b_0=1$ as well. Its products of one-forms are all zero, whereas the triangular calculus has volume form $D(u)D(h)D(d)$.
\end{example}

\begin{example}[The enveloping algebra of $\mathfrak{sl}_2$]\label{sl2-smooth}
The down-up algebra $A(2,-1,-2)$ is differentially smooth. More precisely, its realization as $U(\mathfrak{sl}_2(\Bbbk))$ admits a connected integrable calculus with right module ranks $1,6,6,1$ in degrees $0,1,2,3$, respectively.
\end{example}
\begin{proof}
Let $\mathfrak g=\mathfrak{sl}_2(\Bbbk)$ with basis $e,f,H$ and brackets $[e,f]=H$, $[H,e]=2e$, and $[H,f]=-2f$. Put $U=U(\mathfrak g)$, which is isomorphic to $A(2,-1,-2)$ by Example \ref{sl-two}.

We first construct the algebra of invariant forms. Take two copies of the exterior algebra $\bigwedge\mathfrak g$, with the common degree-zero part $\Bbbk$ and with all products between their positive-degree parts equal to zero. This is the fiber product of the two exterior algebras over their augmentations. Denote the two copies of $e\wedge f\wedge H$ by $v_1,v_2$. Their difference spans a two-sided ideal because both volume elements annihilate the positive-degree part. Let $\mathcal E$ be the quotient by this ideal and write $v$ for their common image. Thus
\begin{align*}
\mathcal E^0&=\Bbbk,&\quad
\mathcal E^1&=\mathfrak g\oplus\mathfrak g,&\quad
\mathcal E^2&=\bigwedge^2\mathfrak g\oplus\bigwedge^2\mathfrak g,&\quad
\mathcal E^3&=\Bbbk v.
\end{align*}
The algebra $\mathcal E$ is graded-commutative and is generated by $\mathcal E^1$. For $j=1,2$, multiplication gives a nondegenerate pairing $\mathcal E^j\times\mathcal E^{3-j}\to\Bbbk v$, since in each copy it is the usual exterior pairing and the cross terms vanish.

The adjoint action of $\mathfrak g$ on each copy of $\mathfrak g$ extends to an action by degree-zero derivations on the exterior algebras. The action fixes $v_1,v_2$ since $\operatorname{tr}(\operatorname{ad}_x)=0$ for each $x\in\mathfrak g$. It therefore descends to $\mathcal E$. Form the graded algebra $\Omega=\mathcal E\rtimes U$, with underlying space $\mathcal E\otimes_{\Bbbk} U$ and commutation rules
\[
x\eta=\eta x+x\cdot\eta,\qquad x\in\mathfrak g,\quad \eta\in\mathcal E,
\]
where $x\cdot\eta$ denotes the induced adjoint action. Explicitly, the coproduct $\Delta(x)=x\otimes 1+1\otimes x$ on $U$ gives the product
\[
(\eta\otimes a)(\xi\otimes b) =\sum_{(a)}\eta\bigl(a_{(1)}\cdot\xi\bigr)\otimes a_{(2)}b, \quad \Delta(a)=\sum_{(a)}a_{(1)}\otimes a_{(2)}.
\]
It is associative by coassociativity and the identity $a\cdot(\eta\xi)=\sum_{(a)}(a_{(1)}\cdot\eta)(a_{(2)}\cdot\xi)$, which follows from the derivation property on the generators of $U$. In particular, $\Omega^0=U$ and $\Omega^j=\mathcal E^j\otimes U$ as a right $U$-module.

Put $\theta=(H,e)\in\mathcal E^1$. For homogeneous $\omega\in\Omega$, define
\[
D(\omega)=\theta\omega-(-1)^{|\omega|}\omega\theta.
\]
This is a graded derivation and $D^2(\omega)=\theta^2\omega-\omega\theta^2=0$, since $\theta^2=0$ in $\mathcal E$. In particular, $D(x)=-x\cdot\theta$ for $x\in\mathfrak g$, so
\[
D(e)=(2e,0),\quad D(f)=(-2f,H),\quad D(H)=(0,-2e).
\]
The $U$-bimodule generated by $D(U)$ is stable under $x\cdot\eta=x\eta-\eta x$ on its invariant forms. Starting from $D(e)$ and $D(H)$ and acting twice by $f$, we obtain all of $\mathfrak g$ in each copy because $[f,e]=-H$ and $[f,-H]=-2f$. Hence $U D(U) U=\Omega^1$. As $\mathcal E$ is generated in degree one, $\Omega$ is generated by $U$ and $D(U)$. The Leibniz rule moves intervening coefficients to the left, so the density condition in Definition \ref{calculus} holds in every degree.

To prove connectedness, filter $U$ by the usual PBW filtration $F_mU$ and set $F_m\Omega^j=\mathcal E^j\otimes F_mU$. The three displayed values of $D$ imply $D(F_mU)\subseteq F_{m-1}\Omega^1$. The commutation rule with invariant forms shows that its leading-symbol map acts on $P\in S(\mathfrak g)$ by
\[
P\longmapsto D(e)\otimes\frac{\partial P}{\partial e}+D(f)\otimes\frac{\partial P}{\partial f}+D(H)\otimes\frac{\partial P}{\partial H}.
\]
The three vectors $D(e),D(f),D(H)$ are linearly independent over $\Bbbk$. Consequently, if $a\in U$ has PBW degree $m>0$, then the leading homogeneous symbol of $a$ has some nonzero partial derivative and $D(a)\ne0$. We conclude that $\ker(D|_U)=\Bbbk1$.

Finally, any $\Bbbk$-basis of $\mathcal E^j$ is also a left $U$-basis of $\Omega^j$. Indeed, the map $U\otimes\mathcal E^j\to\Omega^j$ given by left multiplication is filtered and induces the tensor flip on associated graded spaces. Since the filtrations are exhaustive and bounded below, the map is bijective. The element $v$ commutes with $U$, so $\Omega^3=vU=Uv$ with volume automorphism $\nu_v=\mathrm{id}_U$. Fix $j=1,2$ and bases $\{\eta_i\}$ of $\mathcal E^j$ and $\{\xi_i\}$ of $\mathcal E^{3-j}$ chosen so that $\eta_i\xi_k=\delta_{ik}v$. Given $\sum_i a_i\eta_i\in\Omega^j$ in its left basis, we have
\[
\ell_{\pi_v}^{j}\!\left(\sum_i a_i\eta_i\right)(\xi_k)=a_k.
\]
Since the $\xi_k$ are a right basis, these evaluations show that $\ell_{\pi_v}^{j}$ is bijective. By Proposition \ref{integrability}, the calculus is integrable. The PBW theorem gives $\operatorname{GKdim}U=3$, which is the dimension of the calculus.
\end{proof}

The same conclusion holds for $A(2,-1,\gamma)$ for every $\gamma\in\Bbbk$. For $\gamma=0$ this is Example \ref{heisenberg-smooth}. For $\gamma\ne0$, replacing $d$ by $-2\gamma^{-1}d$ and keeping $u$ changes the parameter to $-2$, as follows directly from \eqref{du-first} and \eqref{du-second}.

The calculus in Example \ref{sl2-smooth} has six independent invariant one-forms. The obstruction for semisimple Lie algebras in \cite[Section 4]{MartiniPostKersten1995} assumes that the differentials of a Lie basis themselves form the invariant one-form basis. Over $\mathbb{C}$, the four-dimensional calculus of Batista and Majid \cite{BatistaMajid2003} on the angular momentum algebra supplies a further comparison. Our calculus has top degree three, as required by Definition \ref{smoothness}.

\subsection{Generalized down-up algebras}

\begin{example}
Let $r,s\in\Bbbk^*$. For $f=0$ and $\varepsilon=0$, Theorem \ref{gdu-constant} gives the scalar exterior calculus on the quantum affine  in Example \ref{gdu-examples} (i). The same formulas with $T(h)=rh-\varepsilon$ give a calculus for every $\varepsilon$. For instance, $L(0,1,2,1)$ is differentially smooth with $T(h)=h-1$.
\end{example}

\begin{example}[Weyl and quantum Weyl relations]
For $c\ne0$, the relation $du-ud=-c$ gives
\[
L(c,1,1,0)\cong A_1(\Bbbk)\otimes_{\Bbbk}\Bbbk[h],
\]
where $A_1(\Bbbk)$ denotes the first Weyl algebra. Theorem \ref{gdu-constant} gives its calculus and also applies to $L(1,2,3,1)$, in which $du-3ud=-1$ and $T(h)=2h-1$. In the latter case the own-generator assignments are $\nu_u(u)=3^{-1}u$ and $\nu_d(d)=3d$. The factors $[n]_3$ and $[n]_{3^{-1}}$ occurring in the proof are nonzero for all $n\geq1$.
\end{example}

\begin{example}\label{gdu-monomial}
Suppose $r\ne1$, put $h_0=\varepsilon/(r-1)$, and let
\[
f(x) = c(x-h_0)^m,\quad c\ne0,\quad m\geq2.
\]
Then $f(T(h))=r^m f(h)$. If $q\in\Bbbk^*$ satisfies $q^m=r^m/s$, take $G(h)=h_0+q(h-h_0)$ and $\lambda=r^m/s$. The maps $G$ and $T$ commute, and $f(G(h))=\lambda f(h)$. Theorem \ref{triangular-gdu} applies if $q$ and $\lambda$ are each $1$ or of infinite multiplicative order. In particular, $s=r^m$ allows $q=\lambda=1$, so Corollary \ref{cor:gdu-eigen} applies directly. A case with both mixed commutations present is $L((X-1)^2,2,1,1)$, for which $h_0=1$, $q=2$ and $\lambda=4$. Since $f$ has degree two, Proposition \ref{gdu-diagonal-obstruction} excludes the diagonal construction on these three generators.
\end{example}

\begin{example}
For $r=s=1$ and $\varepsilon=0$, Corollary \ref{cor:gdu-eigen} applies to every $f\in\Bbbk[X]$. Here $h$ is central, the forms satisfy the ordinary exterior relations, and the only nontrivial commutator of a generator with a basis form is
\[
dD(u)-D(u)d=-D(h)f'(h).
\]
This gives, for instance, $L(X^2+1,1,1,0)$ without an eigenvalue restriction on the coefficients of the polynomial.
\end{example}

\section{Scope of the constructions and future work}\label{Futurework}

Propositions \ref{gdu-diagonal-obstruction} and \ref{du-diagonal-obstruction} give obstructions under specified bimodule hypotheses. The following argument explains the restriction in the ordinary down-up case for a larger class of affine commutation rules.

\begin{proposition}\label{du-affine-obstruction}
Put $x_1=d$ and $x_2=u$ in $A=A(\alpha,\beta,\gamma)$. Suppose a connected differential calculus over $A$ satisfies
\[
x_iD(x_j)=D(x_j)(a_{ij}x_i+b_{ij}),\quad a_{ij}\in\Bbbk^*,\quad b_{ij}\in\Bbbk,\quad i,j\in\{1,2\}.
\]
Then $\Omega^m A=0$ for all $m\geq3$. Thus these affine rules cannot give a calculus witnessing differential smoothness of $A$.
\end{proposition}
\begin{proof}
Put $e_i=D(x_i)$. If $a_{ii}=-1$, then $D(x_i^2-b_{ii}x_i)=0$. The element $x_i^2-b_{ii}x_i$ is nonscalar by Proposition \ref{du-pbw}, contrary to connectedness. Applying $D$ to the rule for $x_i e_i$ now gives $(1+a_{ii})e_i^2=0$, so $e_i^2=0$. Similarly the rule for $x_i e_j$ gives $e_i e_j=-a_{ij}e_j e_i$ for $i\ne j$. The Leibniz rule expresses $D(A)$ in the bimodule generated by $e_1,e_2$. The affine rules move all intermediate coefficients to the right. By density, every higher form is consequently a sum of words in $e_1,e_2$ with right coefficients. A word of length at least three has a repeated factor, which can be moved next to its other occurrence and gives zero. Since $\operatorname{GKdim} A=3$, the calculus cannot witness differential smoothness.
\end{proof}

Examples \ref{heisenberg-smooth} and \ref{sl2-smooth} give counterexamples to \cite[Theorem 3.3]{RubianoReyes2026BQAThree}. To identify the two rows in \cite[Table 1]{RubianoReyes2026BQAThree}, take $(x_1,x_2,x_3)=(d,u,du-ud)$ in both cases. Here $q_1=q_2=q_3=1$. Then $[x_2,x_1]=-x_3$, so $c=-1$ in the notation of \cite[Proposition 2.11]{RubianoReyes2026BQAThree}. In the Heisenberg case, $[x_3,x_1]=[x_3,x_2]=0$. In the $\mathfrak{sl}_2$ case, $[x_3,x_1]=2x_1$ and $[x_3,x_2]=-2x_2$. Thus both rows satisfy the hypothesis $c\ne0$ of the earlier theorem, and the connected integrable three-dimensional calculi constructed here show that both entries should read \emph{differentially smooth} in characteristic zero.

The proof of the earlier theorem expresses $D(x_3)$ in the bimodule generated by $D(x_1)$ and $D(x_2)$. The step from this bimodule generation statement to $\Omega^3=0$ requires control of the coefficients between one-forms. In Example \ref{heisenberg-smooth}, the mixed relation $dD(u)-D(u)d=-D(h)$ exhibits $D(h)$ as a bimodule combination while $D(u),D(h),D(d)$ remain independent over the algebra on the right. Proposition \ref{du-affine-obstruction} gives a vanishing result under the additional affine diagonal commutations stated there.

There remain parameter ranges for which the constructions of Section \ref{Differentialandintegralcalculusdownup} do not decide differential smoothness. If $r/s$ has finite order $m>1$ in Theorem \ref{TheoremDU}, the displayed triangular calculus satisfies $D(u^m)=0$. The monomial $u^m$ is nonscalar, so this calculus is disconnected. A separate question is whether another calculus gives smoothness in such cases, for instance for $A(0,1,0)$, which also appears in \cite[Remark 4.7]{BrzezinskiLomp2018}. Other cases of interest are $\beta=0$, $rs=0$, and generalized down-up algebras for which $f$ does not satisfy \eqref{eigen-gdu}. The generalized Weyl presentation \eqref{gdu-gwa} and the skew derivations of \cite{AlmulhemBrzezinski2019} provide further structures to examine in these cases.

\section{Declarations}

The authors have no conflict of interest to disclose.

\end{document}